\documentclass[11pt]{article}

\usepackage{amsmath, amssymb, amsthm, enumerate}
\usepackage{graphicx, caption, subcaption}
\usepackage{xcolor}
\usepackage{tikz-cd} 
\usepackage[pagebackref]{hyperref}  

\usepackage{cleveref}

\usepackage[margin=1.25in]{geometry}

\usepackage{titlesec}
\theoremstyle{plain}

\theoremstyle{definition}

\newcommand{\C}{\mathbb{C}}

\newcommand{\R}{\mathbb{R}}

\newcommand{\D}{\mathbb{D}}

\newcommand{\M}{\mathcal{M}}

\newcommand{\T}{\mathcal{T}}

\renewcommand{\H}{\mathcal{H}}
\renewcommand{\P}{\mathbb{P}}

\newcommand{\Mod}{\operatorname{Mod}}

\newcounter{bencomments}

\AtBeginDocument{%
   \def\MR#1{}
 }

\title{Compactifications of strata of differentials}
\author{Benjamin Dozier\thanks{Department of Mathematics, Cornell University, \href{mailto:benjamin.dozier@cornell.edu}{\nolinkurl{benjamin.dozier@cornell.edu}}.}} 
\begin{document}
\maketitle

\begin{abstract}
In this informal expository note, we quickly introduce and survey compactifications of strata of holomorphic 1-forms on Riemann surfaces, i.e. spaces of translation surfaces.   In the last decade, several of these have been constructed, studied, and successfully applied to problems.  We discuss relations between their definitions and properties, focusing on the different notions of convergence from a flat geometric perspective. 
\end{abstract}

\section{Motivation}
\label{sec:motivation}


The study of holomorphic differentials on Riemann surfaces has been a vibrant area of research for many years.  Differentials have persistently played a major role in the algebro-geometric study of curves and in Teichm\"uller theory.  The cotangent space of Teichm\"uller space at a Riemann surface $S$ is naturally identified with holomorphic quadratic differentials on $S$.   A somewhat newer aspect is the connection to billiard dynamics.  A translation surface is a Riemann surface equipped with a holomorphic $1$-form; one source of these is unfoldings of rational billiard tables. There is a $GL_2(\R)$ action on the space of translation surfaces, which is central to understanding dynamics on individual translation surfaces as well as Teichm\"uller geodesic flow.
These spaces, which are stratified according to the multiplicity of zeros of the differential, have rich structure, and admit various interesting types of geometry.

A fundamental property of these strata is that they are \emph{non-compact}.  This makes the theory interesting, but also can be a major obstacle; for instance, it means that recurrence properties of various flows are not automatic.  Given a non-compact space, it is often fruitful to study its \emph{compactifications}.  Here are several reasons for this:
\begin{enumerate}
\item Shedding light on the nature of the non-compactness of the original space.  For instance, it might contain information about the number of ends of the original space. 
\item  \label{item:pd} Access to results for compact spaces.  This occurs for the basic example of affine space $\C^n$, which can be compactified by projective space $\C\P^n$.  In $\C\P^n$, one has Bezout's theorem on intersections as well as tools such as Poincar\'e duality for compact manifolds. 
These algebraic and topological motivations apply to strata of differentials.  
\item  The recursive structure exhibited by the boundary.  Each boundary we will discuss itself parametrizes Riemann surfaces and differentials -- the compactification is \emph{modular}.  The surfaces on the boundary are more degenerate, hence simpler.  This makes possible certain inductive arguments, provided one can achieve suitable degenerations. 
\end{enumerate}

Because of these motivations, and others, there has been interest in understanding the compactifications that strata admit.  This topic has been studied in some form for many decades.  But in recent years, the pace of developments has increased.  The compactifications have played central roles in several recent advances in the theory of translation surfaces.  Our knowledge about these compactifications is beginning to reach a mature stage.   The resulting spaces, though enjoying beautiful properties, are somewhat complicated even to define.  The goal of this short article is to convey the main ideas behind several of them as succinctly as possible, and to focus on the similarities and differences between them.  We assume the reader is already interested in strata of differentials and has some basic knowledge of them.  

\paragraph{Plan of the article.}

We will first recall in \Cref{sec:deligne-mumford} the natural Deligne-Mumford compactification of the moduli space of Riemann surfaces.  Next we move on to strata of differentials, discussing the sources of non-compactness in \Cref{sec:sources}.   We then focus on four recently constructed compactifications of strata of differentials.   These naturally fit into a hierarchy in terms of how much information is remembered about diverging sequences.  We discuss the compactifications in order of increasing amount of information (which also roughly corresponds to increasing complexity of definition of the space).  We conclude in \Cref{sec:overview} with a summarizing table and chart.  



\section{Acknowledgements}
\label{sec:ack}

I would like to thank Francois Greer, Samuel Grushevsky, and Alex Wright for helpful conversations.  I also thank Barak Weiss for useful comments on an earlier draft.  And I thank the anonymous referees for helpful suggestions.  My research is supported in part by NSF Grant DMS-2247244 and the Simons Foundation.  

\section{Deligne-Mumford compactification}
\label{sec:deligne-mumford}

\begin{figure}
\centering
  \includegraphics[scale=0.36]{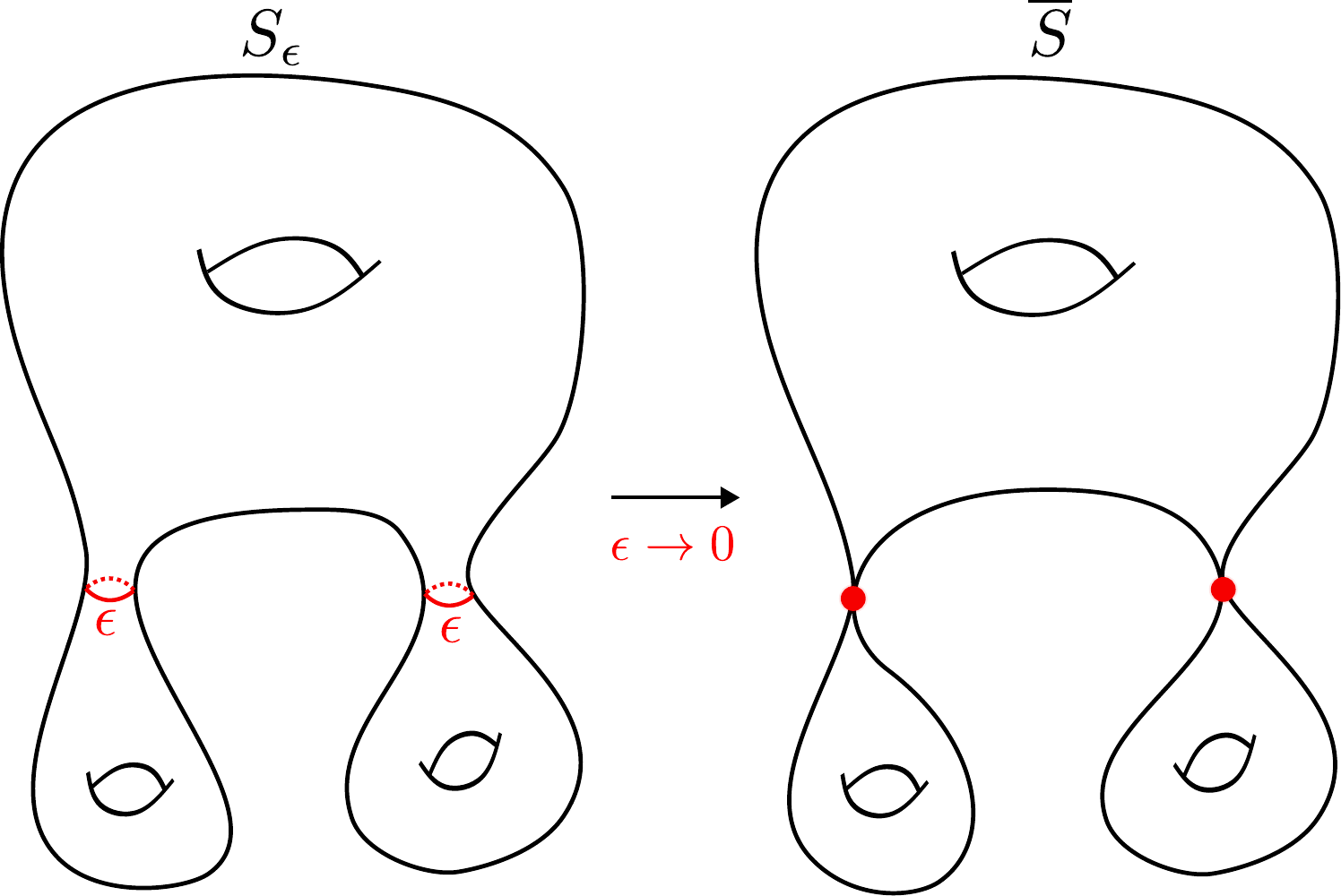}
  \caption{Riemann surfaces in $\M_3$ with red curves getting pinched converge to a nodal Riemann surface in $\overline\M_3$.  The limit surface lies in a piece of $\partial \overline \M_3$ parametrized by $\M_{1,2}\times \M_{1,1} \times \M_{1,1}$ (actually, it is the quotient of this by the group with two elements permuting the last two factors; this is not a serious issue, though it also comes up for compactifications of strata of differentials, and we will ignore it in the rest of the paper).}
  \label{fig:mgbar}
\end{figure}

We begin with a brief description of the Deligne-Mumford \cite{dm69} compactification $\overline \M_g$ of $\M_g$, the space of genus $g$ Riemann surfaces (without any differentials).   This turns out to be really the ``one'' compactification of $\M_g$.  Most of the compactifications of strata of differentials are based in some way on Deligne-Mumford.  It can be defined from the perspective of algebraic geometry, hyperbolic geometry, or complex analysis.  We will discuss the hyperbolic perspective; the complex analytic one is similar, while the algebro-geometric one is rather different and involves geometric invariant theory.  The basic idea behind the compactification is that certain curves on the surface can get ``short''; in the limit these curves are pinched and the Riemann surface develops a node.  See \Cref{fig:mgbar}.

To make this slightly more precise, one introduces the \emph{augmented Teichm\"uller space} $\widehat \T_g$, which is a non-compact bordification of Teichm\"uller space $\T_g$ (the orbifold universal cover of $\M_g$, parametrizing Riemann surfaces endowed with a marking up to homotopy).  From any pants decomposition of $S_g$ (a topological surface of genus $g$), we get a system of global Fenchel-Nielsen coordinates.  For each cuff of the pants, there is a (hyperbolic) length parameter taking values in $\R_{>0}$ and a twist parameter taking values in $\R$.  Then $\widehat \T_g$ is defined by adding points to $\T_g$ for which some cuff lengths in a Fenchel-Nielsen chart are zero; this space has a natural topology.  The new points do not contain any information about the twist parameter for cuffs that get pinched.  For more details, see \cite{abikoff77} and references therein.

The space $\widehat \T_g$ is not particularly nice; for instance it is not locally compact, since any neighborhood of a boundary point contains surfaces with any twist parameter in $\R$ for a pinched cuff curve.  However, the action of the mapping class group $\Mod(S_g)$ on $\T_g$ extends nicely to $\widehat \T_g$.  The Deligne-Mumford compactification $\overline \M_g$ is defined to be the quotient, and it is a very nice space.

\paragraph{Properties.}  The Deligne-Mumford space has the structure of a smooth compact complex orbifold.  In fact, via the algebro-geometric perspective, it's a smooth projective variety.
\footnote{though because special Riemann surfaces admit automorphisms, in many cases it should be thought of as a \emph{stack}}

To explain the smoothness of $\M_g$ at the boundary, note that we can think of length and twist as like magnitude and angle in polar coordinates for $\C$ (the twist for $\T_g$ took values in $\R$, but after quotienting by the action of the mapping class group, in particular the Dehn twist about the pinched curve, we are left with a value in the circle).   
Then adjoining points with zero length parameter is the same as adding in the origin of the polar coordinate system.

The compactness can be seen by combining (i) the existence of constants $B_g$ (the Bers constants) such that any surface in $\M_g$ has a pants decomposition with all cuffs of hyperbolic length at most $B_g$, with (ii) the fact that in each genus there are only finitely many types of pants decompositions, up to the action of the mapping class group.     

The boundary $\partial \M_g := \overline \M_g - \M_g$ is essentially a union of certain products of moduli spaces of lower complexity Riemann surfaces.  These Riemann surfaces have punctures corresponding to the pinched curves; for this reason it is better to start out by generalizing the problem to compactifying $\M_{g,n}$, the moduli space of genus $g$ Riemann surfaces with $n$ punctures/marked points.   The boundary  $\partial \M_g$ is not quite a smooth subvariety of $\M_g$; it has several irreducible components, which can intersect each other, and some of the components have self-intersections.  However, $\partial \M_g$ stills sits rather nicely in $\M_g$; it is a \emph{normal crossing divisor}.  

\paragraph{Applications.}  The construction of the projective algebraic variety $\overline \M_g$ shows that $\M_g$ is a quasi-projective variety, and thus it has solid footing in the world of algebraic geometry.  

Using $\overline \M_g$, it was shown in the landmark paper \cite{mh82} that $\M_g$ is a variety of \emph{general type} for $g\ge 24$;  informally, this means that there is no way of parameterizing the set of Riemann surfaces of such genus in a nice way.    This was considered a big surprise at the time.

The compactification has also been used to great effect in Brill-Noether theory, which studies maps from algebraic curves to projective spaces of various dimensions (in e.g. \cite{gh80}).  One often can show that curves with a certain property form a Zariski open subset of $\overline \M_g$, but it can be challenging to show that such a set is actually non-empty.  Boundary points of $\overline \M_g$, which are nodal curves, can be helpful in this regard, since, when sufficiently degenerate, they become combinatorial (and often tractable) to analyze.  

Deligne-Mumford space also has various applications in Teichm\"uller geometry.  For instance, in a recent paper \cite{ds21} of the author and Sapir, it's used to show that no proper algebraic subvariety of $\M_g$ (in particular, a projection of a stratum of differentials) can be coarsely dense with respect to the Teichm\"uller metric.

\section{Sources of non-compactness of strata}
\label{sec:sources}
We now equip our surfaces with differentials.  We define the \emph{stratum} $\H=\H(\kappa_1,\ldots,\kappa_m)$ to be set of holomorphic differentials on a compact Riemann surface with $m$ labeled points where the differential vanishes to orders $\kappa_1,\ldots,\kappa_m$ (and is non-vanishing everywhere else).   These are the basic objects we will compactify.  We note here that not all strata are connected, but each has at most $3$ connected components, and the classification of these components is well-understood \cite{kz23}.  

We will generalize the above definition somewhat, allowing the $\kappa_i$ to take any integer value (such strata will naturally appear in our boundaries).  An entry $\kappa_i=0$ indicates a \emph{marked point}.  An entry $\kappa_i<0$ indicates we are allowing a \emph{meromorphic differential}, and it has a pole of order $\kappa_i$ at the $i$th marked point (and is holomorphic away from the marked points).  


The non-compactness of $\H$ has three sources discussed in the next subsections.  The fact that these are the only such sources follows from the thick-thin decomposition for flat surfaces (see \cite[Section 4.3]{minsky} and \cite{rafi}).

\subsection{$\R_{>0}$ action}
\label{sec:r-action}
The most obvious source of non-compactness of $\H$ is the $\R_{>0}$ action by scaling the differential.  It is generally not useful to add a limit point to a sequence such as $X, 2X, 3X,\ldots$, where $X$ is some fixed element of $\H$.  Instead we will consider compactifications of $\P\H$, the quotient of $\H$ by the $\C^*$ action (in what follows $\P$ will always denote a quotient by $\C^*$). Sometimes it is better to quotient by $\R_{>0}$, or, equivalently, to consider only unit-area surfaces.  Quotienting by $\R_{>0}$ retains more information than quotienting by $\C^*$, but the disadvantage is that the resulting space loses its complex analytic structure.   

\subsection{High modulus flat cylinders}
\label{sec:high-mod-cyl}
\begin{figure}
\centering
  \includegraphics[scale=0.66]{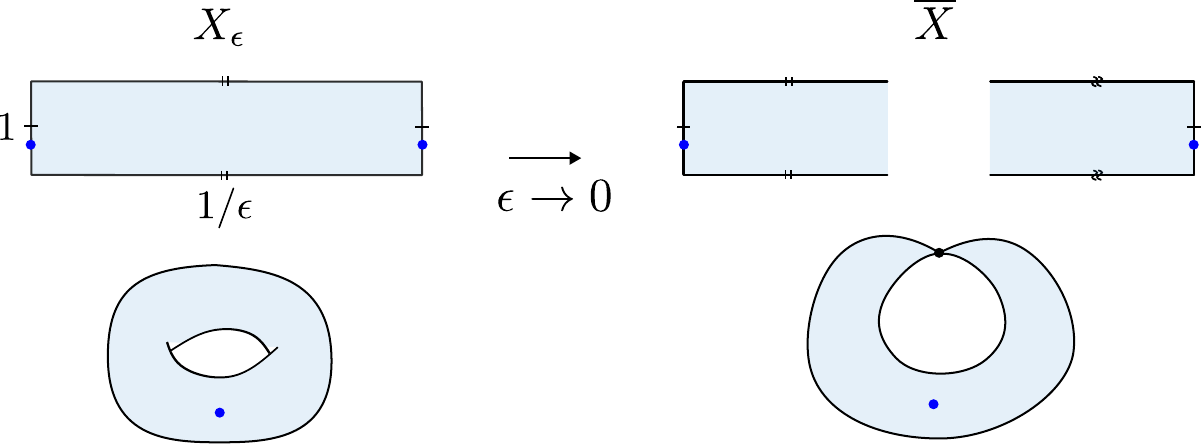}
  \caption{Surfaces in $\H(0)$ converging to a meromorphic differential with a pair of simple poles.  The top shows convergence of flat structures, while the bottom depicts convergence of the underlying Riemann surfaces.}
  \label{fig:degen_cyl}
\end{figure}

A sequence of surfaces in $\H$ with Euclidean cylinders of aspect ratio going to infinity diverges. In fact, the sequence of underlying Riemann surfaces already diverges, since the cylinders give conformal annuli of high modulus.   This phenomenon already occurs in $\mathcal{H}(0)$, the space of genus $1$ flat surfaces with a single marked point; see \Cref{fig:degen_cyl}. 

All of the compactifications we will consider deal with this source of non-compactness in the same way.  The flat limit of the cylinders is taken to be two half-infinite Euclidean cylinders.  The underlying Riemann surfaces converge to a singular Riemann surface with a node corresponding to the infinite ends. 

We can also consider the nodal Riemann surface with the node punctured out, and then fill in with two marked points, giving a smooth compact Riemann surface.  The limit of the holomorphic $1$-forms develops \emph{simple poles} at the marked points.  One can see how a simple pole corresponds to a half-infinite cylinder by pulling back the $1$-form $dz$ (which descends to a $1$-form on the Euclidean half-infinite cylinder) to the punctured disc $\D^*$ via the map $z\mapsto \log z$.  The new $1$-form is $d\log z = \frac{1}{z}dz$, which has a simple pole at $0$.  This point should be thought of as corresponding to the point ``at the end'' of the half-infinite cylinder.

\subsection{Small subsurfaces}
\label{sec:small-subsurf}
\begin{figure}[]
\begin{center}
  \includegraphics[scale=0.75]{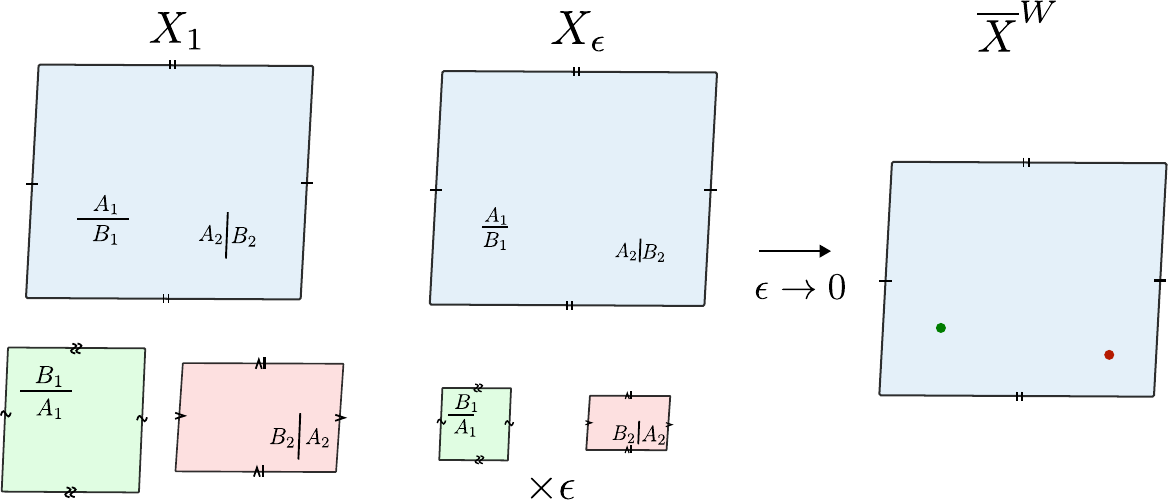}
  \caption{Convergence of a family of translation surfaces in the WYSIWYG $\overline \H^W$.  The surfaces $X_\epsilon$ are in the genus $3$ stratum $\H(1,1,1,1)$, while the limit lies in the genus $1$ stratum $\H(0,0)$.  Since the red and green subsurfaces get very small, they do not appear in the limit.      The green and red points on $\overline X^{W}$ record the location where the green and red subsurfaces ``disappeared''.  }
  \label{fig:degen_W}
\end{center}
\end{figure}

One can construct translation surfaces by gluing along slits, and thus it is possible for a whole subsurface to be very small.  The running example we will focus on is the family given in the middle of \Cref{fig:degen_W}, in which the whole green and red subsurfaces are getting small as $\epsilon\to 0$.  There is no analogous phenomenon in the world of hyperbolic surfaces, since the collar lemma forbids two short hyperbolic closed geodesics from intersecting.

The rest of this article is devoted to describing how the different compactifications deal with small subsurfaces, and the differences in properties that result.  We define them in order of increasing amount of information remembered about limit points.

\section{What You See Is What You Get}
\label{sec:wysiwyg}

The WYSIWYG compactification is in some sense the most straightforward, though it has the disadvantage that the total space is not particularly nice.  The notion of convergence we will describe was introduced in \cite{mcmullen}, and the WYSIWYG space was extensively studied in \cite{mw17}.


The reader should begin by studying \Cref{fig:degen_W} to see what convergence in WYSIWYG is meant to capture.  

The way that the WYSIWYG deals with small subsurfaces is to simply forget them.  This is the sense of the word ``See'' in What You See Is What You Get; since the small parts become microscopic, one does not see them in the limit.  Given polygon representations $X, X_1,X_2,\ldots$ of translation surfaces, perhaps in different strata (and we also allow disconnected surfaces, and marked points), there is a natural notion of geometric convergence $X_n\to X$: edges of the polygons should converge in the limit.  The WYSIWYG space $\overline \H^W$ is then defined by attaching to $\H$ any surface with marked points (from any stratum) that can be obtained as such a limit of surfaces in $\H$.  Additionally, one attaches surfaces with infinite cylinders (i.e. compact Riemann surfaces with a meromorphic differential all of whose poles have order at most $1$), as described in \Cref{sec:high-mod-cyl}. This gives $\overline \H^W$ as a set, and then one endows it with a natural topology.  \footnote{A formal definition of the convergence above goes as follows: there should be maps $f_n: X-U_n\to X_n$, diffeomorphisms onto their images, where $U_n$ are a decreasing family of open neighborhoods of the marked points of $X$, such that (i) the pullbacks of the differentials from $X_n$ converge to the differential on $X$, and (ii) the injectivity radius of points in $X_n - f_n(X-U_n)$ tends to zero.}

Note that the limit of a sequence can be disconnected.  For instance, start with the surface on the left in \Cref{fig:degen_W}, but now shrink the upper blue part, keeping the red and green subsurfaces fixed size (and shrinking both slits).  The limit has two components, corresponding to the red and green subsurfaces; it lies in a piece of the boundary parametrized by $\H(0)\times \H(0)$.

The boundary of $\overline \H^W$ is thus a union of certain strata of translation surfaces (possibly disconnected and with marked points).  


\paragraph{Properties.}
The quotient $\P \overline \H^W$ is a compactification of $\P \H$.  Informally, this is because all three sources of non-compactness described in \Cref{sec:sources} have been handled.  With some effort, this can be made rigorous; alternatively, a quick proof can be given once the Hodge bundle compactification is defined -- see discussion of WYSIWYG in \Cref{sec:hodge-bundle}.   

Since $GL_2(\R)$ acts on all the objects in the definition, there is a natural continuous extension of its action from $\H$ to all of $\overline \H^W$.  

The WYSIWYG is built out of pieces from various strata.  Although each stratum has nice properties (in particular, each is a quasi-projective algebraic variety), the way that they are glued together turns out to make the resulting space not so nice.  In particular, there is a natural map from an algebraic variety, the projective Hodge bundle compactification $\P \overline \H^{HB}$ (see next section) to $\P \overline \H^W$ ,  but there is no way to give $\P \overline \H^W$ the structure of an algebraic variety in such a way that this map is a morphism of varieties \cite[Theorem 1.1]{cw21}.

Unlike for many of the other spaces we will study, there is not a natural map $\overline \H^W\to \overline \M_g$.  This is because, the WYSIWYG really forgets all information about parts of the surface that are getting vary small, while Deligne-Mumford remembers the conformal structure there.  

\paragraph{Applications.}
The WYSIWYG compactification has been used to great effect in the study of $GL_2(\R)$-orbit closures in $\H$, a central topic in Teichm\"uller dynamics.  Since the $GL_2(\R)$ action extends to the WYSIWYG boundary, the closure in the boundary can be understood using the same type of tools that apply in $\H$, in particular the results of Eskin-Mirzakhani \cite{em18} and Eskin-Mirzakhani-Mohammadi \cite{emm15}, which give that each $GL_2(\R)$-orbit closure is an \emph{affine invariant manifold} $\mathcal N$.

Provided that one can show certain degenerations exist within $\mathcal N$, an inductive classification strategy using $\overline \H^W$ becomes possible.  Using this, Mirzakhani-Wright \cite{mw18} classified full rank affine invariant manifolds.  This was then generalized in a series of papers, culminating with \cite{aw23}, which classifies high rank affine invariant manifolds.

\section{Hodge bundle compactification}
\label{sec:hodge-bundle}

\begin{figure}[]
\begin{center}
  \includegraphics[scale=0.75]{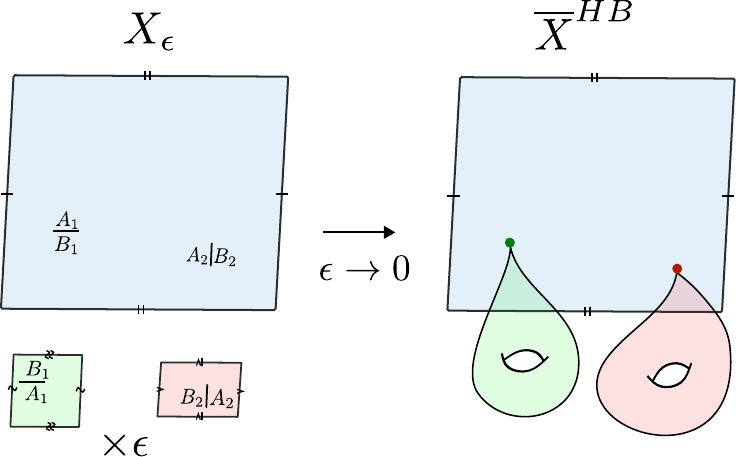}
  \caption{Convergence in the Hodge bundle $\overline \H^{HB}$.  On $\overline X^{HB}$ the blue piece is a component on which the differentials converge to a non-zero differential. On the green and red components, the differentials tends to $0$, and only the Riemann surface structure here is remembered.  The limit $\overline X^{HB}$ lies in a piece of the boundary parametrized by $\H(0,0) \times \M_1 \times \M_1$.}
  \label{fig:degen_HB}
\end{center}
\end{figure}

The Hodge bundle compactification $\P \overline \H^{HB}$ records all the data of the WYSIWYG, while not completely ignoring information about small subsurfaces.  Instead, the limit of the underlying Riemann surfaces of such vanishing subsurfaces is remembered, while any extra information from the differential is forgotten.   See \Cref{fig:degen_HB}.

The above essentially already defines $\overline \H^{HB}$, but to give a more formal description, we begin by recalling the definitions of some related spaces.  The \emph{Hodge bundle} $\Omega \M_g$ is the holomorphic vector bundle over $\M_g$ whose fiber over a Riemann surface $S$ is the vector space of holomorphic $1$-forms on $S$.  This bundle has a natural extension to a holomorphic vector bundle $\Omega \overline \M_g$ over $\overline \M_g$, which we also call the Hodge bundle.  These bundles have played a major role in algebraic geometry for many decades; for instance they can be used to prove the projectivity of $\overline \M_{g,n}$ \cite{cornalba93}, and to exhibit interesting cohomology classes on moduli space \cite{mumford83}. The fiber over a nodal Riemann surface $\overline S$ consists of certain meromorphic $1$-forms on $\overline S$ with at worst simple poles. The condition is that the differential is holomorphic, except possibly at the nodes, where the differential is allowed to have simple poles; in this case, there must be a simple pole on both sides of the node, and the residues must be opposite.  In terms of flat geometry, these correspond to pairs of half-infinite cylinders, as in \Cref{fig:degen_cyl}.

Note that the differential is allowed to vanish on some (or even all) components of $\overline S$.   If a sequence of elements in $\Omega \M_g$ converges to such a limit in $\Omega \overline \M_g$,  then there are subsurfaces on which the differential is tending to $0$.  The boundary point only remembers the limiting conformal structure here -- no further flat information is recorded.   

The space $\overline \H^{HB}$ can then be succinctly defined as the closure of the stratum $\H$ in $\Omega \overline \M_g$.  All of these spaces admit a $\C^*$ action by scaling of differentials, and $\P \overline \H^{HB}$ is just the quotient by this action.  Since $\Omega \overline \M_g$ is a vector bundle over a compact space, $\P \Omega \overline \M_g$ is compact, and hence so is $\P \overline \H^{HB}$.

The boundary of $\overline \H^{HB}$ is a union of certain products of strata of translation surfaces (with marked points) and moduli spaces (with marked points) of Riemann surfaces.

\paragraph{Relation to WYSIWYG.} We could have also defined the WYSIWYG using the Hodge bundle: $\overline \H^W$ is the quotient of $ \overline \H^{HB}$ that identifies two boundary surfaces if the result of removing components where the differentials vanish, and filling in the resulting punctures with marked points, yields the same translation surfaces with marked points.  From this perspective, it is easy to see that $\P \overline \H^W$ is in fact compact.

\paragraph{Properties.}

The space $\P \overline \H^{HB}$ is a projective algebraic variety, since it is the closure of a quasi-projective variety inside the projective variety $\P \Omega \overline \M_g$.
\footnote{The closure is taken in the Euclidean topology, but for a quasi-projective variety this is the same as the Zariski closure.}
However $\overline \H^{HB}$ is not a smooth variety (see discussion of end of \Cref{sec:ivc} -- that also applies to the Hodge bundle).  



\paragraph{Applications.} Besides it's natural importance in algebraic geometry, the Hodge bundle compactification (or rather, a variant of it for quadratic differentials) has also been used in \cite{bdr24, aw24} to study boundaries of subvarieties of moduli space that are totally geodesic with respect to the Teichm\"uller metric.

\section{Incidence variety compactification}
\label{sec:ivc}

As we move up in our hierarchy of compactifications, we would like to remember more flat geometric information about the shape of small subsurfaces.  Since the differentials are converging to zero, the idea is to rescale them so they converge to a non-zero limit.  One imagines having a set of magnifying glasses of various magnifications to look at various parts of the surface as they are degenerating.  In this sense one remembers more than what one ``sees'' with just  the naked eye in the WYSIWYG.   

For the IVC compactification, we start with the Hodge bundle compactification, and then additionally remember a differential, up to complex scaling, for each component of the limiting Riemann surface on which the differential vanishes.  This space was defined and studied by Bainbridge-Chen-Gendron-Grushevsky-M\"oller in \cite{bcggm18}. 

\begin{figure}[]
\begin{center}
  \includegraphics[scale=0.65]{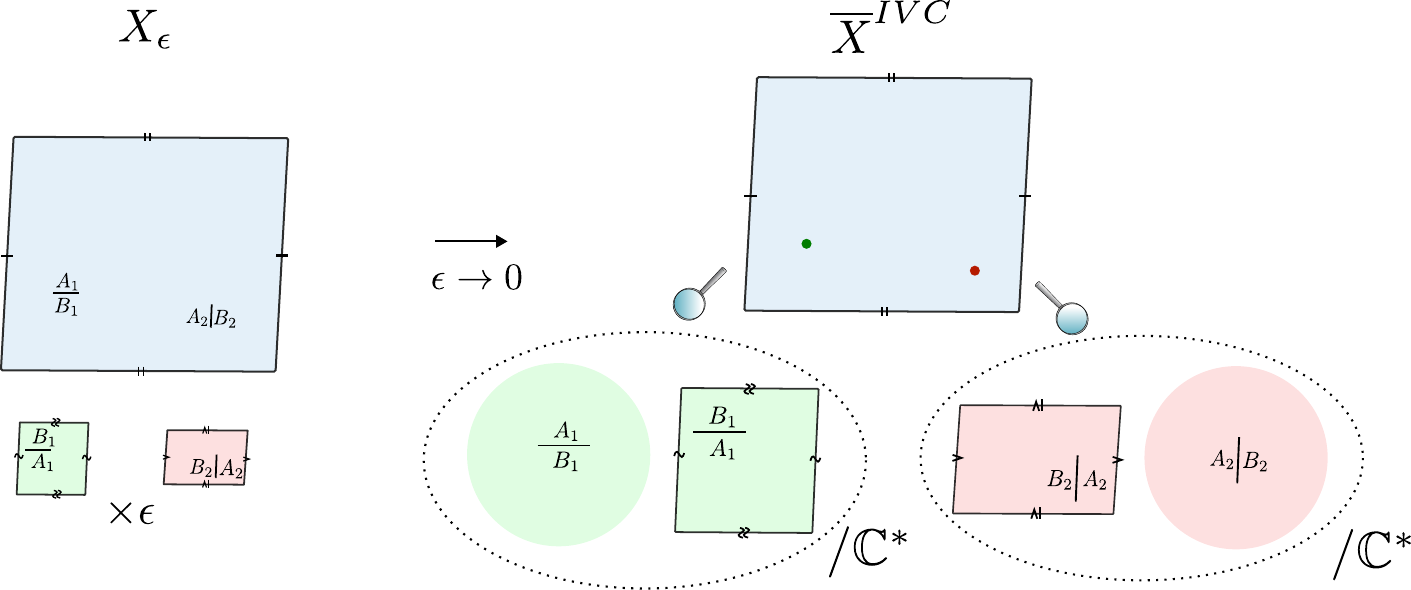}
  \caption{Convergence in the incidence variety compactification $\overline \H^{IVC}$.   On the upper, blue component of $\overline X^{IVC}$, it remembers a holomorphic differential with marked points.  For each of the green and red surfaces, it remembers a Riemann surface with 2 marked points recording the location of the zeros of the degenerating differentials.  On each, we get a unique \emph{meromorphic} differential, up to complex scaling, with a pole at the node, and zeros at the marked points. Flat pictures of these differentials are shown on the bottom right; each represents a flat torus glued into a copy of $\C$ (a portion of which is drawn as a disc) along a slit.    Each can be thought of as a limit when the component is rescaled to converge to something not identically zero.  The limit $\overline X^{IVC}$ does not remember the relative sizes of, say, the slits on the green and red subsurfaces.
    The piece of the boundary that $\overline X^{IVC}$ lies in is parametrized by $\H(0,0) \times \M_{1,3} \times \M_{1,3}$ (for each $\M_{1,3}$ factor, two of the special points come from zeros of the differential, while the third comes from the pole). } 
  
  \label{fig:degen_IVC}
\end{center}
\end{figure}

In \Cref{fig:degen_IVC} we illustrate how our running example of degenerating family converges in the IVC.

To construct the IVC formally, we start first by defining a slightly simpler compactification $\P \overline \H^{DM}$, based on the Deligne-Mumford compactification $\overline \M_{g,n}$ with marked points.  We leverage the basic complex analytic fact that a holomorphic $1$-form on a smooth Riemann surface is determined, up to complex scale, by the locations of its zeros.  This implies that for a stratum $\H$ with $n$ zeros, there is an embedding $\P \H \to \M_{g,n}$ given by $[(S,\omega)] \mapsto (S,z_1,\ldots,z_n)$, where $z_i\in X$ is the location of the $i$th zero of $\omega$. \footnote{To define this, we should really start with a finite cover of the stratum $\H$ where the zeros are \emph{labeled}, which allows us to sensibly speak of the ``$i$th zero''; this is a minor issue.}  We then define $\overline\H^{DM}$ to be the closure of $\H$ in $\overline \M_{g,n}$ via this embedding.  A point on the boundary is a nodal Riemann surface with marked points, and this combination of data contains information about the shapes of a sequence of translation surfaces converging to the boundary point. 

The IVC is constructed by combining the information of $\P \overline \H^{DM}$ with that of the Hodge bundle compactification $\overline \H^{HB}$.  We consider the embedding $\H \to \Omega \M_{g,n}$ given by $(S,\omega) \mapsto (S, z_1,\ldots,z_n, \omega)$, where again $z_i\in X$ is the location of the $i$th zero of $\omega$.  Then the IVC $\overline \H^{IVC}$ is defined to be the closure of $\H$ in $\Omega \overline \M_{g,n}$ via this embedding.

\paragraph{Interpretation of boundary points.}

Points in $\partial \overline \H^{IVC}$ remember information about differentials that are going to zero on some component; however this information is not directly encoded in a differential, but rather through the limiting locations of the zeros.  In particular, it is not immediately apparent which points in $\Omega \overline \M_{g,n}$ actually lie in the closure of $\H$.

The main result of \cite{bcggm18} is a characterization of the closure in terms of differentials on the components of the limit.  A point $\Omega \overline \M_{g,n}$ lies in the closure of $\H$ if it admits a \emph{twisted differential}, defined as a collection of meromorphic differentials, one on each component of the underlying nodal Riemann surface, satisfying certain explicit conditions.  One condition is that the poles of the differentials only occur at nodes.  The boundary point only remembers this meromorphic differential up to \emph{complex scaling} on each component separately.  In addition to poles of order $1$ arising from high-modulus cylinders, poles of higher order also arise, which give flat surfaces of infinite area -- such a pole represents the larger part of the surface, seen from the perspective of the smaller part.  The model for a pole of order $2$ is the point $\infty$ for the differential $dz$ on the Riemann sphere $\widehat \C$; one computes the order using the coordinate change $\zeta=1/z$.
\footnote{An equivalent description: one pulls back the differential $dz$ on $\widehat \C$ by the map $z\mapsto z^n$, which is an unramified cover away from $0,\infty$.  The resulting differential on $\widehat \C$ has a pole of order $n+1$ at $\infty$. }
The green (resp. red) subsurface on $\overline X^{IVC}$ in \Cref{fig:degen_IVC} has an order $2$ pole, which is the point at infinity in the plane that the green (resp. red) disc lies in.      

Poles of order higher than $2$ also arise.  The flat picture of such a pole looks like multiple copies of $\C$ glued together along ray slits.  Poles that arise often have residues (though this does not occur for the $\overline X^{IVC}$ in our running example).  There are various compatibility conditions on pole orders and residues, notably the Global Residue Condition, needed for a collection of meromorphic differentials to appear in $\partial \overline \H^{IVC}$.

\paragraph{Properties.}  We have defined $\overline \H^{IVC}$ as the closure of $\H$ within an ambient algebraic variety, so it's also an algebraic variety, and in fact $\P \overline \H^{IVC}$ is a projective algebraic variety.  The term \emph{incidence variety} is used because a point in the image of $\H$ has the property that the marked points of the Riemann surface coincide with zeros of the differential, an incidence condition.

However, $\overline \H^{IVC}$ is not a smooth variety; it has fairly serious \emph{singularities}.
It is perhaps not so easy to see this directly.  Morally, it is because the IVC is still ``forgetting too much information'', specifically the information of the relative scales of differentials below the top level.   Once one has the multi-scale compactification, which is smooth, one can use the forgetful map to the IVC to prove non-smoothness of the latter.  The map contracts a locus of a certain dimension which is inconsistent with the target being smooth.   See \cite[Example 14.1]{bcggm19}.   (This is in some sense analogous to the relation between algebraicity of the Hodge bundle and the WYSIWYG -- there is a map from the former to the latter, but it contracts in such a way that the latter cannot be algebraic.)

\paragraph{Applications.}

The IVC has been fruitfully used to compute Masur-Veech volumes of strata of differentials \cite{sauvaget18, cmsz20}.   In these papers, the volumes are related to intersection numbers of certain cohomology classes over the IVC (so this is an instance of \Cref{item:pd} discussed in \Cref{sec:motivation}).  The authors also give various formulas for Siegel-Veech constants (see \Cref{sec:other-bdy} for a discussion of these and the general approach to computing them in terms of volumes).  

In a different application of a more algebro-geometric flavor, \cite{mullane17} uses strata of differentials to construct an infinite family of natural effective divisors on $\overline \M_{g,n}$ and then uses the description of the IVC to show that these are extremal and rigid.  

\section{Multi-scale compactification}
\label{sec:multiscale}

Degenerating families of differentials often have a natural set of scales associated with different parts of the surface.  For instance, in our running example, there is a scale (size around $1$) associated with the blue subsurface, and a smaller scale (size around $\epsilon$) associated with the smaller green and red subsurfaces.   We would like a compactification that (i)  faithfully remembers the set of scales and their ordering, and (ii) also allows us to compare sizes/lengths of objects on subsurfaces of the same scale (e.g. on the green and red subsurfaces in our running example).  The IVC is a step towards these goals, but does not fully achieve either -- for instance the limit example in \Cref{fig:degen_IVC} fails (ii).  This is because the IVC relies on recovering differentials, up to complex scaling, from the locations of their zeros; it does not directly record meromorphic differentials.

\begin{figure}
\begin{center}
  \includegraphics[scale=0.68]{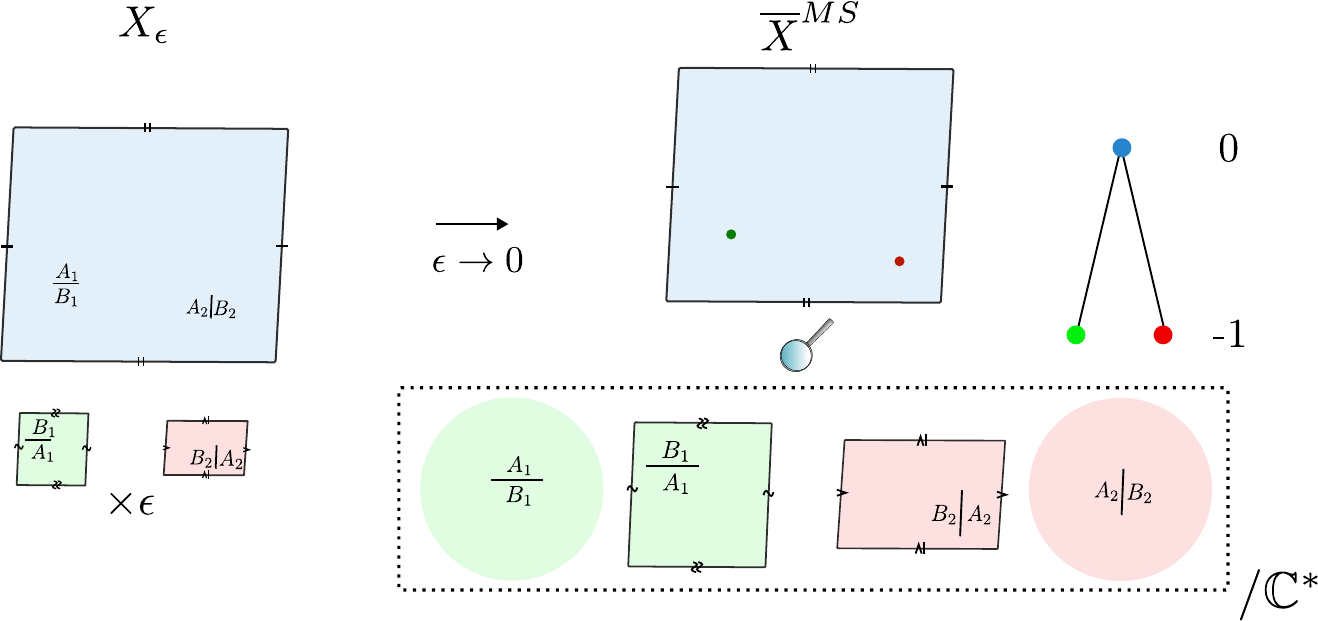}
  \caption{Convergence in the multi-scale $\overline \H^{MS}$.  On the upper blue component of $\overline X^{MS}$, there's a holomorphic differential with marked points. On the union of the lower green and red components, there's a meromorphic differential, up to complex scaling.  The limit $\overline X^{MS}$ remembers the relative sizes of, say, the slits on the green and red subsurfaces, since they lie on the same (disconnected) surface that is projectivized.  The level graph is shown on the right (it is convenient to take the top level to be $0$).  
The limit $\overline X^{MS}$ lies in a piece of the boundary parametrized by $\H(0,0) \times \P[\H(1,1,-2)\times \H(1,1,-2)]$.}
 \label{fig:degen_MS}
\end{center}
\end{figure}

\begin{figure}
\begin{center}
  \includegraphics[scale=0.68]{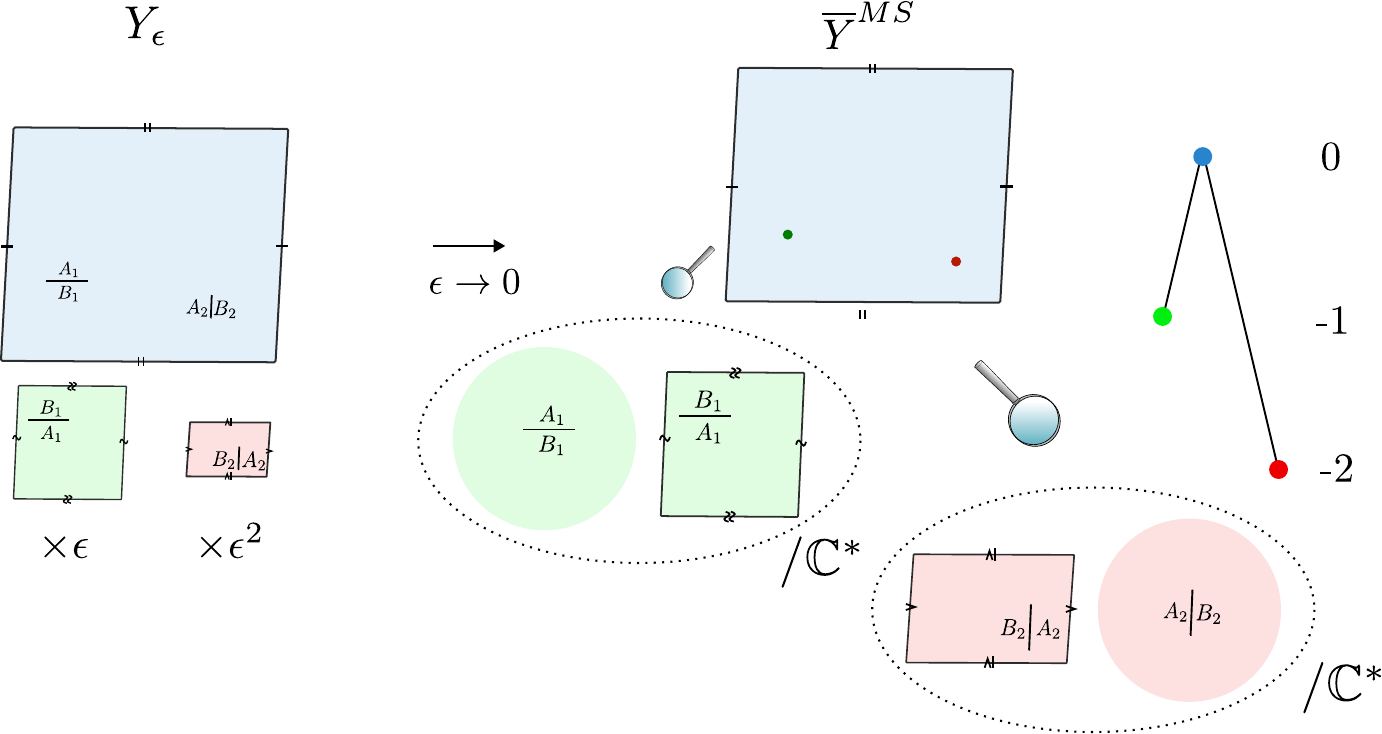}
  \caption{Convergence of a slightly different family $Y_\epsilon$ in $\overline \H^{MS}$; the red subsurface is now getting small much faster than the green, since the red is of size $\epsilon^2$, while the green is of size $\epsilon$.   The level graph has three levels, which encodes that blue was much larger than green which was much larger than red.   The limit data consists of three differentials, one for each of blue, green, and red subsurfaces. The red and green are projectivized separately, since they are at different levels. The limit $\overline Y^{MS}$ lies in a piece of the boundary parametrized by $\H(0,0) \times \P\H(1,1,-2)\times \P\H(1,1,-2)$.
    The limit of $Y_\epsilon$ in the WYSIWYG, Hodge bundle, and IVC is the same as the limit of $X_\epsilon$ in each of these; only in multi-scale do we get different limits.}
  \label{fig:degen_MS_3level}
\end{center}
\end{figure}

The idea of the multi-scale compactification, defined and studied in \cite{bcggm19}, is to fix these deficiencies of the IVC by directly remembering meromorphic differentials, as well as the ordered set of scales.  See \Cref{fig:degen_MS} for how this works in our running example, and \Cref{fig:degen_MS_3level} for a different family that has three different scales.  Unlike the previous spaces we've discussed, $\overline \H^{MS}$ does not admit a succinct description (that we know of); its construction takes a great deal of effort.  Higher order poles and residues both play important roles, as in the IVC.

\paragraph{Boundary points.} The boundary points in the multi-scale compactification are given by (equivalence classes of) \emph{multi-scale differentials}.  These are described by some combinatorial data, a \emph{level graph}, together with the continuous data of a twisted differential. 

The level graph is, as a graph, just the dual graph of the underlying stable Riemann surface (the same underlying Riemann surface as in IVC).  That is, the vertices correspond to components, and for each node joining a pair of components, there is an edge connecting the corresponding vertices.  Additionally, there is a discrete set of levels, indexed by integers, and each vertex of the graph appears at some level; higher levels correspond to larger subsurfaces.

The twisted differential must be \emph{compatible} with the level graph; a key part of this is that for each node, the zero/pole orders of the differentials respect the level structure of the graph in the sense that the differential on the lower component has a pole, and the sum of the orders of vanishing on the two sides is equal to $-2$.  The additional properties that must be satisfied, notably the Global Residue Condition, are generally similar to IVC.    There is an additional subtlety that only arises with poles of order at least $3$ called \emph{prong matching}; it has to do with different combinatorial choices for how to glue in such a differential into a higher component.

Finally, the equivalence relation on multi-scale differentials is given by orbits of the action of the \emph{level rotation torus}.  The group is $(\C^*)^N$, where $N$ is the number of levels, and it acts by rescaling each level (and also takes into account the prong-matching).

Having specified the set of boundary points, to fully define $\overline \H^{MS}$ one then has to give a complex analytic structure to the whole space, which takes a considerable amount of work.  

\paragraph{Properties.}

The payoff for the hard work of constructing $\overline \H^{MS}$ is that we get a very nice space.  In particular, encoding so much information in limit points results in a space that is \emph{smooth}, as a complex orbifold.  The boundary $\partial \overline \H^{MS}$ is a normal crossing divisor.  These properties were proved in the original paper \cite{bcggm19}.   The first construction was complex-analytic and did not yield that $\P \overline \H^{MS}$ is a projective algebraic variety, but this has since been shown in two different ways \cite{cmm22,cghms22}.  

\paragraph{Applications.}

The multi-scale compactification has already played a central role in computations of Chern classes and Euler characteristics of strata of differentials \cite{cmz22}; the smoothness of the space is crucial, since it allows one to invoke Poincar\'e duality.  

This compactification has also proved useful in the study of affine invariant manifolds.  In work of the author \cite{dozier23}, it is used to control the measure of the set of surfaces in an affine invariant manifold that have multiple short saddle connections.  Here the complex-analytic structure and compactness play key roles, as does the property that, since $\P \overline \H^{MS}$ remembers so much information, neighborhoods of boundary points are quite small, which restricts the geometry of the surfaces that lie in them.   Using similar techniques, Bonnafoux \cite{bonnafoux} shows that the Siegel-Veech transform of any affine measure is in $L^2$.  In \cite{bdg22}, $\overline \H^{MS}$ is used to give a new proof and generalization of Wright's Cylinder Deformation Theorem, as well as to better understand the equations that cut out an affine invariant manifold.  

Along these lines, it is hoped that $\overline \H^{MS}$ will be helpful in the classification of affine invariant manifolds via inductive arguments on the boundary.  The WYSIWYG has already been used in this way, as discussed above.  Since the multi-scale remembers more information, different avenues of attack may be available.  

In \cite{benirschke}, double ramification loci, which are of great interest in algebraic and symplectic geometry, are interpreted via exact differentials as linear submanifolds of strata, and their closure is studied via the multi-scale compactification.   Very recently in \cite{cms23}, $\overline \H^{MS}$ is used to investigate Chern classes of affine invariant manifolds.  

\section{Other boundary notions}
\label{sec:other-bdy}

\paragraph{Principal boundary.}

The \emph{principal boundary} of a stratum was constructed by Eskin-Masur-Zorich \cite{emz03} well before the compactifications described above had been defined.  They attach to $\H$ limits of families of surfaces obtained by shrinking configurations of parallel saddle connections.  Their space does not give a compactification; their purpose was rather to compute Siegel-Veech constants, which are averages of geometric counts over strata.  Via the Siegel-Veech integral formula, such a constant can be reinterpreted in terms of the volume of the locus where certain saddle connections are $\epsilon$-short.  This locus can be thought of as a tubular neighborhood of a degenerate stratum where those saddle connections have been contracted to zero length.  Hence the desired volume is the volume of the degenerate stratum multiplied by a factor involving $\epsilon$.  These degenerate strata form the principal boundary.  

The principal boundary has been interpreted in \cite{cc19} as a certain nice subset of $\overline \H^{IVC}$ (and it can also be embedded into the multi-scale $\overline \H^{MS}$; see \cite[Section 4]{lee2023}).  The boundary surfaces here all have at most two levels, and simple poles can only appear on the bottom level.  The short configuration of saddle connections ends up on the lower level.  Each component $C$ of the complement of the configuration on the prelimit surfaces corresponds to a pair of simple poles on the bottom level when $C$ is a flat cylinder, and to a component of the top level otherwise.   

\paragraph{Real multi-scale compactification.}
This is a variant of the multi-scale compactification $\overline \H^{MS}$ defined above, where the quotients by $\C^*$ are replaced by quotients by $\R^*$.  It can also be constructed as a level-wise real oriented blow-up of $\overline \H^{MS}$ (see \cite[Section 12]{bcggm19}).  So a limit point remembers more information; it remembers the angles of periods, even at lower levels.

One advantage of the real multi-scale space is that it admits a continuous extension of the $GL(2,\R)$ action on $\H$.  The space is also compact.  This space (or rather a variant of it for quadratic differentials) and these properties are used in \cite{bdr24} to study boundaries of totally geodesic subvarieties of moduli space.  The downside of this space is that it not a complex algebraic variety; rather it is a an orbifold with \emph{corners}.  

\section{Overview}
\label{sec:overview}

We summarize the properties of and relations between the compactifications in the following table and diagram.  

\begin{table}[h]
\begin{tabular}{|l|l|l|l|}
\hline
  & $GL_2(\R)$ action extends & \begin{tabular}{@{}c@{}} Projectivization is \\ a projective variety \end{tabular} & Smooth orbifold  \\ \hline
  WYSIWYG $\overline \H^{W}$ &  Yes &  No &  No \\ 
  Hodge Bundle $\overline \H^{HB}$ & No & Yes  & No\\ 
  IVC $\overline \H^{IVC}$ & No  &  Yes & No  \\ 
  Multi-scale $\overline \H^{MS}$ & No  & Yes  & Yes \\ \hline
\end{tabular}
\caption{Properties of the four compactifications} 
\end{table}

\begin{figure}[h]
  \centering
\[\begin{tikzcd}
	{\P\overline \H^{MS}} & {\P\overline \H^{IVC}} & {\P\overline \H^{HB}} & {\P\overline \H^{W}} \\
	&& {\P\overline \H^{DM}} & {\overline \M_g}
	\arrow[from=1-1, to=1-2]
	\arrow[from=1-2, to=2-3]
	\arrow[from=1-2, to=1-3]
	\arrow[from=1-3, to=1-4]
	\arrow[from=1-3, to=2-4]
	\arrow[from=2-3, to=2-4]
      \end{tikzcd}\]
    \caption{Natural forgetful maps between the various compactifications fit into a commutative diagram.}
    \label{fig:relations}
\end{figure}

{\footnotesize
\bibliographystyle{amsalpha}
  \bibliography{sources}
  }

  
\end{document}